\documentclass[conference]{IEEEtran}
\IEEEoverridecommandlockouts
\usepackage{cite}
\usepackage{amsmath,amssymb,amsfonts, amsthm}
\usepackage{algorithm}
\usepackage{algpseudocode}
\usepackage{graphicx}
\usepackage{textcomp}
\usepackage{xcolor}
\usepackage{bm}
\usepackage{tikz}
\usepackage{tikzscale}
\usepackage{enumitem}
\usepackage{hyperref}

\newcommand*{\CopyCounter}[2]{%
  \expandafter\def\csname c@#2\endcsname{\csname c@#1\endcsname}%
  \expandafter\def\csname p@#2\endcsname{\csname p@#1\endcsname}%
  \expandafter\def\csname the#2\endcsname{\csname the#1\endcsname}}

\newcounter{ModelProblem}
\CopyCounter{Theorem}{Problem}
\CopyCounter{Theorem}{Proposition}
\CopyCounter{Theorem}{Property}
\CopyCounter{Theorem}{Claim}
\CopyCounter{Theorem}{Lemma}
\CopyCounter{Theorem}{Corollary}
\CopyCounter{Theorem}{Conjecture}
\CopyCounter{Theorem}{Definition}
\CopyCounter{Theorem}{Example}
\CopyCounter{Theorem}{Remark}
\CopyCounter{Theorem}{Question}
\CopyCounter{Theorem}{Condition}
\CopyCounter{Theorem}{Criterion}
\CopyCounter{Theorem}{Observation}
\theoremstyle{plain}

\newtheorem{problem}[Problem]{Problem}

\newtheorem{modelProblem}[ModelProblem]{Model Problem}

\newlist{lambdalist}{itemize}{1}
\setlist[lambdalist]{label=\textbf{Q}:}
\newcommand\itemC{\item[\textbf{$C$}:]}
\newcommand\itemB{\item[\textbf{$B$}:]}

\begin{document}

\title{\LARGE \bf
Optimality of Motion Camouflage Under Escape Uncertainty
}

\author{Mallory E. Gaspard$^{1,2}$
\thanks{*Supported by NSF DMS award 1645643.}
\thanks{$^{1}$Email: meg375@cornell.edu}
\thanks{$^{2}$Center for Applied Mathematics, Cornell University, Ithaca, NY 14853}
\thanks{\tiny{© 2024 IEEE. Personal use of this material is permitted. Permission from IEEE must be
obtained for all other uses, in any current or future media, including
reprinting/republishing this material for advertising or promotional purposes, creating new
collective works, for resale or redistribution to servers or lists, or reuse of any copyrighted
component of this work in other works.}}
}

\maketitle

\begin{abstract}
This letter proposes a novel continuous-time dynamic programming framework to determine when it is optimal for a pursuer to use MC amidst uncertainty in the evader's escape attempt time. We motivate this framework through the model problem of an energy-optimizing male hover fly pursuing a female hover fly for mating. The time at which the female fly initiates an escape is modeled to occur as the result of a non-homogeneous Poisson point process with a biologically informed rate function, and we obtain and solve two Hamilton-Jacobi-Bellman (HJB) PDEs which encode the pursuer’s optimal trajectories. Our numerical experiments and statistics illustrate when it is optimal to use MC pursuit tactics amidst uncertainty and how MC optimality is affected by certain properties of the evader's sensing abilities.
\end{abstract}

\section{Introduction}
\label{section:intro}

\begin{figure}[h]
\centering
\includegraphics[width = 0.2\textwidth]{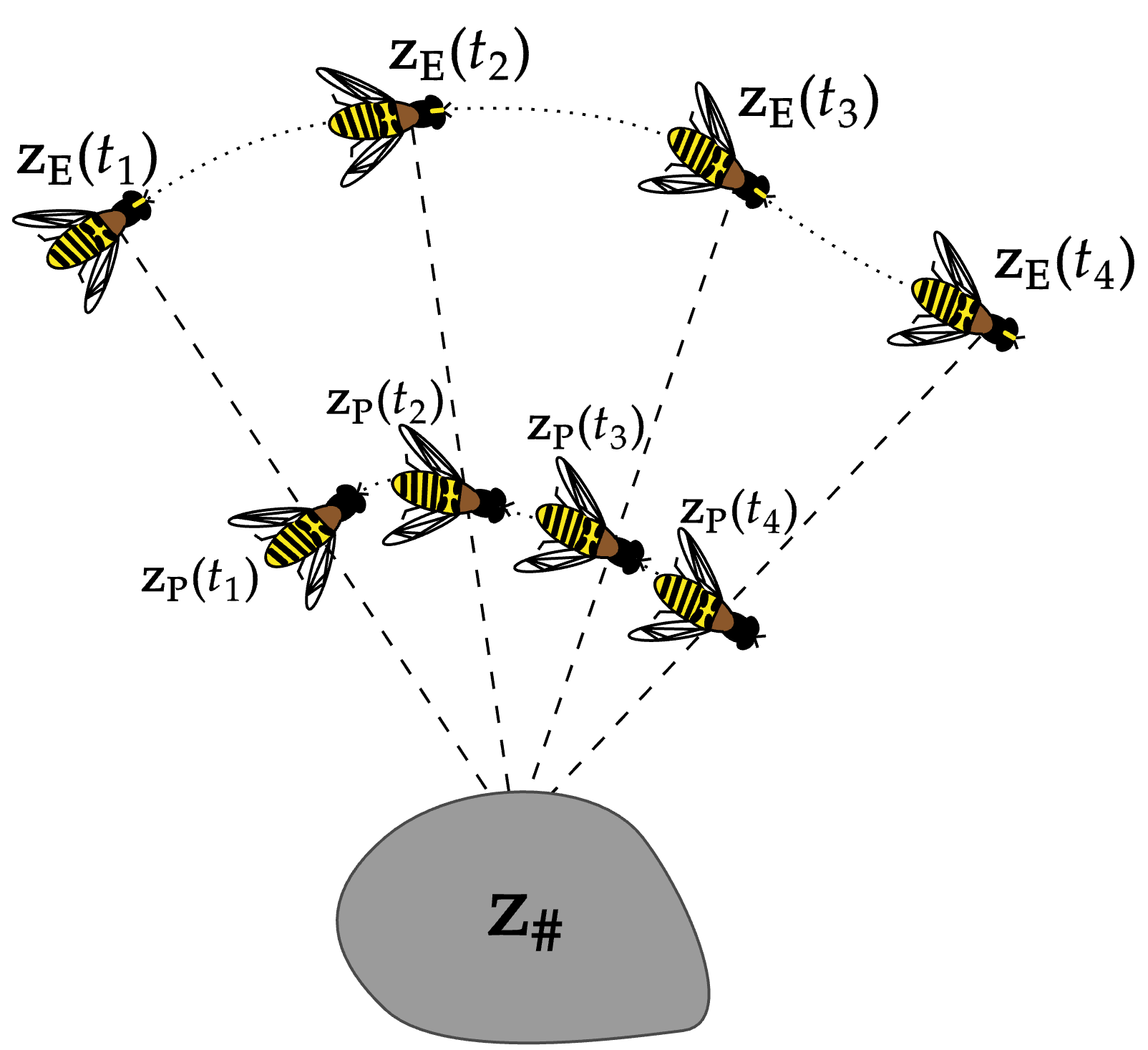} 
\caption{Illustration of stationary point motion camouflage with respect to the boulder at $\boldsymbol{z}_{\#}$. The evader and pursuer at four snapshots in time are identified by $\boldsymbol{z}_E(\cdot)$ and $\boldsymbol{z}_P(\cdot)$ respectively. 
\label{fig:stat_mc}}
\end{figure}

Camouflaging is often used in nature as a concealment mechanism by pursuers hoping to sneak up on a potential evader. Although organisms across the animal kingdom utilize a variety of such tactics, the active motion camouflage (MC) displayed by male hover flies during mating \cite{collett1975visual} and dragonflies settling territorial disputes \cite{mizutani2003motion} is of particular interest in both biological and engineering applications since it is induced by the pursuer's \emph{movement} rather than by a set physical characteristic such as fur patterns or body texture.
A pursuer engaging in stationary-point MC chooses their trajectory in attempt to trick the evader's visual system into believing that the pursuer is a part of their perceived background optical flow \cite{cuthill2019camouflage}. Figure \ref{fig:stat_mc} illustrates the MC motion constraints geometrically. When the pursuer's trajectory stays on the line between the evader and the inanimate object at the evader's focal point ($\boldsymbol{z}_{\#}$), the evader can only detect changes in the pursuer's size and not any relative motion. 
\\
\hspace{0.4cm}Modeling MC trajectories in a control-theoretic setting has been of interest to biologists and engineers alike since Srinivasan and Davey first proposed a mathematical description of this behavior in \cite{srinivasan1995strategies}. 
Prior studies in the literature which address these tactics in an optimal control framework primarily aim to determine optimal trajectories along which the pursuer always uses MC  \cite{xu2008subspace, rano2012optimal, justh2006steering, carey2004biologically}.
Such trajectories are generally obtained by optimizing a performance index in an open-loop framework under the MC motion constraints. These approaches do not address the fundamental question of when it may be beneficial for the pursuer to utilize MC, and they also rely on the realism-limiting assumption that the evader will never feel threatened and attempt to escape prior to capture.

The main contribution of this letter is the development of a continuous-time dynamic programming framework that can be used to determine \emph{when} it is optimal for a pursuer to utilize MC amidst uncertainty in the evader's escape attempt time. Our framework results in solving a sequence of Hamilton-Jacobi-Bellman (HJB) equations from which we can recover the pursuer's optimal trajectories. MC-optimality from a given starting location is then determined by identifying portions of the optimal trajectory that coincide with an MC-trajectory.

Although our approach is easily adapted to analyze benefits of MC usage in a variety of biological and vehicular systems (see \cite{strydom2017uas, jianqing2023guidance, tan2024antipredator} for examples), we develop our framework through the model problem of an energy-optimizing male hover fly (pursuer) trying to capture a female hover fly (evader) for mating while managing (time and location dependent) uncertainty in whether or not she will feel threatened and attempt to escape. The evader's escape attempt is influenced by rate function $\lambda$ that controls the rate of switching from an initial stalk phase in which the pursuer may utilize MC, to a direct chase phase in which the pursuer rapidly darts after the escaping evader.
We then recover the pursuer's optimal trajectories by solving the resulting sequence of HJB PDEs involving $\lambda$.

This letter is outlined as follows: Section \ref{section:prob_form} describes the model problem, biological assumptions, and the two pursuit phases. Section \ref{section:energy_opt_trajs_under_uncertainty} presents the main optimal control framework and the HJB PDEs. 
Section \ref{section:numerics} outlines the numerical methods used to solve the PDEs. Section \ref{section:experiments} shows four examples which showcase the pursuer's balance between energy conservation and trying to remain concealed. It also reports statistics on the percentage of flight spent in MC from thousands of starting points which illustrate how properties of E's visual system affect the optimality of MC tactics. Conclusions and future directions are discussed in Section \ref{section:conclusions}.

\vspace{-0.6cm}
\section{Model Problem Formulation}\label{section:prob_form}
We motivate the development of our optimal control framework by focusing on the following biological model problem: 
\begin{modelProblem}
When is it optimal for an energy-optimizing male hover fly to engage in stationary point motion camouflage while pursuing a female hover fly who may, at some unknown time, feel threatened and attempt to escape?
\end{modelProblem}
\vspace{-0.5cm}
\subsection{Notation and Assumptions:} 
\textit{Notation:} Throughout the paper, we will refer to the male hover fly (pursuer) as ``P" and the female hover fly (evader) as ``E." We let $\boldsymbol{z}_P = [x_P, y_P]^T$ and $\boldsymbol{z}_E = [x_E, y_E]^T$ refer to points in $\mathbb{R}^2$ which specify P and E's initial positions, while $\boldsymbol{z} = [x,y]^T$ refers to a generic point in $\mathbb{R}^2$.
The vectors $\boldsymbol{z}_P(s) = [x_P(s), y_P(s)]^T$ and $\boldsymbol{z}_E(s) = [x_E(s), y_E(s)]^T$ reflect the time-dependent nature of P and E's positions. 

\textit{Biological Assumptions:}
Each fly is represented as a point mass (with $m = 1$) moving around $\mathbb{R}^2$ with no other flies or external forces present. Since male hover flies are known to patrol resource-rich locations for mates \cite{collett1975visual}, we assume that P pursues E while she is en route to feed at the only food source in the region - a flower located at $\boldsymbol{z}_{\star} = (x_{\star}, y_{\star})$. Male hover flies can predict a female's trajectory reasonably well \cite{collett1978hoverflies}, so we suppose that P has complete knowledge of E's route $\boldsymbol{z}_E(t)$ to $\boldsymbol{z}_{\star}$ at the start of the planning process. 

With compound-eyes, hover flies have $\approx 360^{\circ}$ vision, but they can only see objects within a $D \hspace{0.1cm} m$ radius \cite{collett1978hoverflies}. Once P and E are mutually visible, we assume P maintains sight of E for the rest of the pursuit. If P begins planning their trajectory at a location which is not currently visible to E but will be at a later time, he hovers in place until E becomes visible. Then, he engages in a simplified form of the tracking behavior described in \cite{collett1975visual}. Here, P tracks E by flying toward her with speed $F_P$ until arriving at the tracking distance $\overline{D} \hspace{0.1cm} m$ away from her. Once E arrives at $\boldsymbol{z}_{\star}$, P flies directly toward her with speed $F_P$ in attempt to intercept her there, following the ``constant incomer speed" observation described in \cite{thyselius2018visual}.

Energy Expenditure Quantification:  For simplicity, we focus on quantifying P's ``control energy." We define this quantity as the sum of P's minimum operational energy cost (approximated by constant $W \in \mathbb{R}_+$) and P's instantaneous kinetic energy at time $t$. Mathematically, this is given by 
\begin{equation}
	K(\boldsymbol{v}(t)) = W + \frac{|\boldsymbol{v}(t)|^2}{2}.
\end{equation}
\vspace{-0.5cm}
\subsection{Pursuit Phases and Planning Horizons}
Based on behavioral observations \cite{weihs1984optimal}, the planning process consists of two phases:
\begin{enumerate}
	\item \emph{Stalk Phase}: P is approaching E. E has not yet felt threatened by P and has not initiated an escape attempt. P optimizes for energy loss and may use MC tactics.
	\item \emph{Direct Chase Phase}: E feels threatened by P and immediately attempts to escape. Energy conservation is no longer a priority for P, and a high-speed chase ensues. 
\end{enumerate}
\textit{Stalk Phase Time Horizon \& Quantifying MC Deviation: } The stalk phase begins at $t = 0$ and lasts until $T_s = \min\left\{\hat{T}_1, T_1\right\}$. $\hat{T}_1$ is the uncertain time at which E initiates her escape, and $T_1 = \min\{t \ge 0 \hspace{0.1cm} | \hspace{0.1cm} \boldsymbol{z}_P(t) \in \mathcal{T}\}$, where $\mathcal{T} := \{\boldsymbol{z} \in \mathbb{R}^2 \hspace{0.1cm} | \hspace{0.1cm} |\boldsymbol{z} - \boldsymbol{z}_E(t)| \leq \epsilon\}$.
E's escape attempt can either occur while she travels to $\boldsymbol{z}_{\star}$ or while she's stationary at $\boldsymbol{z}_{\star}$ during feeding. 
The escape attempt time is a non-homogeneous exponentially distributed random variable with pointwise rate function $\lambda$ that depends on P's trajectory at time $t$ and properties of E's visual system. The rate is influenced by P's instantaneous angular displacement from the vector connecting $\boldsymbol{z}_{E}(t)$ to her focal point at $\boldsymbol{z}_{\#}$ and two constants $B > 0$ and $C \ge0$:
\begin{lambdalist}
	\itemB How well E's visual system can resolve angles (acuity).
	\itemC How tolerant E is to P's proximity within the visual field. 
\end{lambdalist}
Larger values of $C$ indicate higher tolerance to P's proximity, while lower values of B indicate that E's visual system has stronger acuity. $B$ can be tuned to reflect species-specific capabilities within the current environment \cite{land1997visual}, while $C$ can be determined from experimental data to reflect E's typical response behavior to other organisms within visual range \cite{thyselius2018visual}. 
At time $t$, the angular displacement between E's focal line to $\boldsymbol{z}_{\#}$ and P's current position $\boldsymbol{z}$ is given by 
\begin{equation}\label{eqn:theta_sharp}
	\theta_{\#}(\boldsymbol{z}, t) := \cos^{-1}\left( \frac{\boldsymbol{r}_{\#}(t) \cdot \boldsymbol{r}(t)}{|\boldsymbol{r}_{\#}(t)| |\boldsymbol{r}(t)|}\right),
\end{equation}
where $\boldsymbol{r}(t) := \boldsymbol{z} -\boldsymbol{z}_E(t)$, and $\boldsymbol{r}_{\#}(t) := \boldsymbol{z}_{\#} - \boldsymbol{z}_E(t)$. We also note that $\theta_{\#}$ is irrelevant when $|\boldsymbol{r}(t)| = 0$ since capture will have already occurred and when $|\boldsymbol{r}_{\#}(t)| = 0$ since E cannot pass through $\boldsymbol{z}_{\#}$. Geometrically, $\theta_{\#}$ is an exact quantification of P's instantaneous MC violation as discussed in \cite{rano2013direct}.

Taking inspiration from a two-dimensional Gaussian, we use $\theta_{\#}$, $B$, and $C$ to define
the pointwise escape rate
\begin{equation}\label{eqn:lambda}
    \lambda(\boldsymbol{z}, t) := 
        \begin{cases}
    A,  & \hspace{-0.2cm} |\boldsymbol{r}(t)| \leq \epsilon \\
    Ae^{-\frac{C(|\boldsymbol{r}(t)| - \epsilon)^2}{B + \theta_{\#}}},  & \hspace{-0.1cm} \epsilon < |\boldsymbol{r}(t)| \leq D \\
        Ae^{-C(|\boldsymbol{r}(t)| - \epsilon)^2},  & \hspace{-0.2cm} |\boldsymbol{r}_{\#}(t)| > D \\
    0, & \hspace{-0.2cm}  |\boldsymbol{r}(t)| > D
    \end{cases}
\end{equation}
where $A$ is a positive constant set to reflect the overall strength of E's visual system. 
The probability that E does not attempt an escape within the time interval $[t_1, t_2]$ is
\begin{equation}\label{eqn:prob_of_no_spook}
p = e^{-\int^{t_2}_{t_1} \lambda(\boldsymbol{z}(s), s) ds}.
\end{equation}
We note that when $t \ge T_{\star}$, the $t$-dependence in $\lambda$ can be suppressed since E's location is fixed at $\boldsymbol{z}_{\star}$.

\begin{figure}[h]
\centering
\includegraphics[width = 0.2\textwidth]{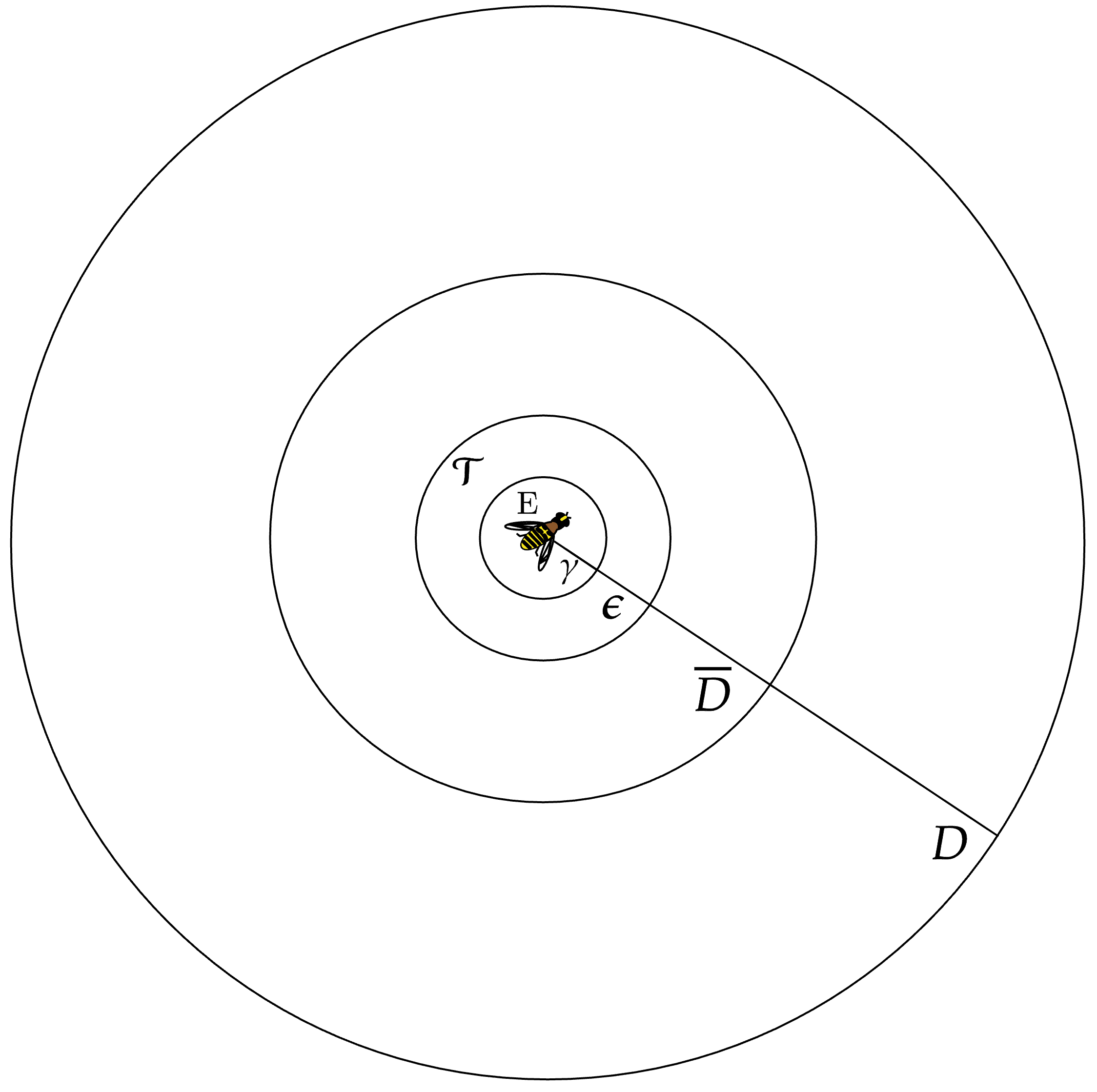} 
\caption{Diagram of all critical distances relative to E.  E cannot see anything beyond $D$ meters. The direct chase phase deterministically begins when P enters $\mathcal{T}$, and E is considered ``captured" when P is within $\gamma$ meters.
\label{fig:mc_dists}}
\end{figure}

\textit{Direct Chase Phase Time Horizon:} At $T_s$, the direct chase begins. Since P usually accelerates to his chase speed at $\approx 500 \hspace{0.1cm} cm/s^2$ \cite{collett1975visual}, we make the simplifying approximation that both P and E instantaneously accelerate to their respective maximum speeds ($G_P$ and $G_E$) at the start of the chase and that $G_P > G_E$. The direct chase phase ends with E's capture at $T_2 = \inf\{t \ge T_s \hspace{0.1cm} | \hspace{0.1cm} |\boldsymbol{z}_P(t) - \boldsymbol{z}_E(t)| \leq \gamma \}$. A visualization of the direct chase region $\mathcal{T}$, and the distances $\gamma, \epsilon, \overline{D}$, and $D$ relative to each other is displayed in Figure \ref{fig:mc_dists}.

\section{Optimal Control Framework}\label{section:energy_opt_trajs_under_uncertainty}
We now present the full optimal control framework that is used to address the Model Problem introduced in Section \ref{section:prob_form}. We begin with an overview of the differential game that models the direct chase phase interaction and then transition into a discussion of the randomly terminated optimal control problem characterizing the stalk phase interaction. 
\subsection{Direct Chase Phase Differential Game}
In the deterministic direct chase phase ($t \ge T_s$), P aims to prevent E's escape. This is modeled by the following problem:
\begin{problem}\label{prob:prob_dc}
	Solve a zero-sum, two-player differential game for the time it takes P to capture E starting from their respective locations, and compute P's direct chase energy expenditure. 
\end{problem}
P freely chases E throughout the plane with dynamics
\begin{equation}\label{eqn:direct_chase_dynamics}
\dot{\boldsymbol{z}}_P(t) = G_P\boldsymbol{a}(t)\hspace{0.2cm} \text{and} \hspace{0.2cm} \dot{\boldsymbol{z}}_E(t) = G_E\boldsymbol{b}(t), \hspace{0.2cm} t \in [T_s, T_2],
\end{equation}
where $\boldsymbol{a} : \mathbb{R} \rightarrow S^1$ and $\boldsymbol{b} : \mathbb{R} \rightarrow S^1$ are measurable control functions specifying P and E's respective directions of motion. This is the simplest form of a pursuit-evasion differential game as discussed in the opening chapters of \cite{isaacs1999differential}. E's capture is guaranteed by our assumption that $G_P > G_E$, and P's optimal direction of motion is $\boldsymbol{a}_* = \frac{\boldsymbol{z}_E(T_s) - \boldsymbol{z}_P(T_s)}{|\boldsymbol{z}_E(T_s) - \boldsymbol{z}_P(T_s)|}$. For P starting from $\boldsymbol{z}$ at $T_s$ and E starting from the position on her known route to $\boldsymbol{z}_{\star}$ at which she felt threatened ($\boldsymbol{z}_E(T_s)$), the time to E's capture $q(\boldsymbol{z}, T_s)$ is known analytically. P's energy expenditure during that timeframe is simply
\begin{equation}
	\delta(\boldsymbol{z}, T_s) = q(\boldsymbol{z}, T_s)\left(W + \frac{G^2_P}{2}\right).
\end{equation}
\vspace{-0.5cm}
\subsection{Stalk Phase Optimal Control Problem}
During the stalk phase of uncertain duration ($t \in [0, T_s)$), P's course of action is found by solving the following problem:
\begin{problem}
Solve a randomly-terminated optimal control problem to determine P's optimal velocity $\boldsymbol{v}(\cdot)$ which minimizes their expected energy loss until the direct chase begins at $T_s$. When P should utilize MC starting from initial location $\boldsymbol{z}_P$ is then determined by identifying which portions of the resulting optimal trajectory satisfy the MC motion constraints. 
\end{problem}

Since E is stationary on the flower until either she is spooked or P gets within $\epsilon$ cm, we compute P's expected energy loss over two stages $t \in [0, T_{\star})$ and $t \ge T_{\star}$. Our approach to this two-stage control problem falls within the general class of randomly terminated optimal control problems developed in \cite{andrews2014deterministic} and similar to multistage examples discussed in \cite{qi2024iu}. 

\textit{Dynamics:} P can predict E's stalk phase trajectory, so her position $\boldsymbol{z}_E(t)$ is given by a known function $\boldsymbol{f}(t)$, $\boldsymbol{f}: [0, T_{\star}] \rightarrow \mathbb{R}^2$. P moves with the simple dynamics
\begin{equation}\label{eqn:stalk_phase_dynamics}
	\dot{\boldsymbol{z}}_P(t) = \boldsymbol{v}(t)
\end{equation}
where $\boldsymbol{v}(t)$ is restricted to the ball $A_F := \{\boldsymbol{v} \in \mathbb{R}^2 \hspace{0.1cm} | \hspace{0.1cm} 0 \leq |\boldsymbol{v}| \leq F_P\}$ when $t \in [0, T_{\star})$ and the circle  $A_S := \{\boldsymbol{v} \in \mathbb{R}^2 \hspace{0.1cm} | \hspace{0.1cm} |\boldsymbol{v}| = F_P\}$ when $t \ge T_{\star}$. Male hover flies  tend to travel with slower speeds while stalking in order to appear less aggressive to potential mates \cite{collett1978hoverflies}, so $F_P << G_P$.

\textit{Indefinite Horizon ($t \ge T_{\star}$): } If E arrives at $\boldsymbol{z}_{\star}$ \emph{without} being spooked, the expected energy loss after $T_{\star}$ is characterized by an indefinite horizon problem which terminates at $\hat{T}_{\star} = \min\{\hat{T}_1, T_1\}$. The expected cost is
\begin{equation}
	\mathcal{C}_{S}(\boldsymbol{z}, \boldsymbol{v}(\cdot)) = \mathbb{E}_{\hat{T}_{\star}} \left\{\int^{\hat{T}_{\star}}_{t} K(\boldsymbol{v}(s)) \hspace {0.1cm} ds + \delta(\boldsymbol{z}_P(\hat{T}_{\star}), \hat{T}_{\star}) \right\}.
\end{equation}
where $t \ge T_{\star}$ and $\delta(\boldsymbol{z}_P(\hat{T}_{\star}), \hat{T}_{\star})$ is the energy cost computed in Problem \ref{prob:prob_dc} from $\boldsymbol{z} = \boldsymbol{z}_P(\hat{T}_{\star})$. The value function is
\begin{equation}
	w(\boldsymbol{z}) = \inf_{\boldsymbol{v}(\cdot)} \left\{\mathcal{C}_{S}(\boldsymbol{z}, \boldsymbol{v}(\cdot)) \right\},
\end{equation}
and by standard arguments in \cite{bardi1997optimal}, $w(\boldsymbol{z})$ is a viscosity solution of the stationary (time-independent) HJB PDE
\begin{equation}\label{eqn:hjb_stalk_stationary}
        0 = \min_{\boldsymbol{v} \in A_S}\left\{ K(\boldsymbol{v}) +  \nabla w(\boldsymbol{z}) \cdot \boldsymbol{v}\right\} + \\
         \lambda(\boldsymbol{z})(\delta(\boldsymbol{z}, T_{\star}) - w(\boldsymbol{z}))
\end{equation}
with the boundary condition $w(\boldsymbol{z}) = \delta(\boldsymbol{z}, T_1)$ for all $\boldsymbol{z} \in \mathcal{T}$.

\textit{Finite Horizon ($t \in [0, T_{\star})$):} Prior to E's arrival at $\boldsymbol{z}_{\star}$, the expected energy loss is described by a finite horizon problem which terminates at $\hat{T}_{F} = \min\{\hat{T}_1, T_{\star}, T_1\}$. It is given by 
\begin{equation}\label{eqn:stalk_expected_cost}
\mathcal{C}_{F}(\boldsymbol{z}, t, \boldsymbol{v}(\cdot)) = \mathbb{E}_{\hat{T}_F} \left\{\int^{\hat{T}_F}_{t} K(\boldsymbol{v}(s)) \hspace {0.1cm} ds + \tilde{\delta}(\boldsymbol{z}_P(\hat{T}_F), \hat{T}_F) \right\}.
\end{equation}
where $\tilde{\delta}(\boldsymbol{z}_P(\hat{T}_F), \hat{T}_F) = \delta(\boldsymbol{z}_P(\hat{T}_F), \hat{T}_F)$ when $\hat{T}_F \neq T_{\star}$, and $\tilde{\delta}(\boldsymbol{z}_P(\hat{T}_F), \hat{T}_F) = w(\boldsymbol{z}(T_{\star}))$ when $\hat{T}_F = T_{\star}$.
The value function for $t \in [0, T_{\star})$ is defined as 
\begin{equation}\label{eqn:vf_stalk_phase}
    u(\boldsymbol{z}, t) = \inf_{\boldsymbol{v}(\cdot)}\left\{\mathcal{C}_F(\boldsymbol{z}, t, \boldsymbol{v}(\cdot))\right\},
\end{equation}
and arguments in \cite{bardi1997optimal} show that $u(\boldsymbol{z}, t)$ is a viscosity solution of the following time-dependent HJB PDE,
\begin{multline}\label{eqn:hjb_stalk}
        0 = u_t(\boldsymbol{z},t) + \min_{\boldsymbol{v} \in A_F}\left\{ K(\boldsymbol{v}) + \nabla u(\boldsymbol{z}, t) \cdot \boldsymbol{v} \right\} + \\
         \lambda(\boldsymbol{z}, t)(\tilde{\delta}(\boldsymbol{z}, t) - u(\boldsymbol{z},t))
\end{multline}
with the terminal condition $u(\boldsymbol{z}, T_{\star}) = \tilde{\delta}(\boldsymbol{z}, T_{\star}) = w(\boldsymbol{z})$ and the boundary condition $u(\boldsymbol{z}, T_1)  = \delta(\boldsymbol{z}, T_1)$ for all $\boldsymbol{z} \in \mathcal{T}$.

\textit{Motion Constraints \& Value Function at Hovering Points:} P's commitment to maintaining sight of E as soon as she becomes visible introduces a motion constraint which prevents P from assuming any position which is currently not visible to E or will never be visible to E. Let $\mathcal{D}(t) := \left\{\boldsymbol{z} \in \mathbb{R}^2 \hspace{0.1cm} : \hspace{0.1cm} |\boldsymbol{z} - \boldsymbol{z}_E(s)| \leq D, \hspace{0.1cm} s \in [t, T_{\star}] \right\}$ denote all of the positions which are currently visible to E or may be at a later time in the process. The set of disallowed states at time $t$ is $\mathcal{B}(t) = \mathbb{R}^2 \setminus \mathcal{D}(t)$, and we set $u(\boldsymbol{z}, t) = + \infty, \hspace{0.1cm} \forall \boldsymbol{z} \in \mathcal{B}(t)$. 

If planning begins from a point $\boldsymbol{z} \in \mathcal{D}(t)$ which is \emph{not} within E's visual range at $t$, P will hover in place until E becomes visible at $\tilde{t} > t$, and P begins tracking. This behavior informs an estimate for $u(\boldsymbol{z}, t)$ at such points given by
\begin{equation}\label{eqn:visual_distance_cost}
\mathcal{C}_{H} =  W(\tilde{t} - t) + p\left(u(\hat{\boldsymbol{z}}, \hat{t}) + K(\hat{\boldsymbol{v}})\hat{\tau}\right)+ (1-p) \tilde{\delta}(\boldsymbol{z}, \tilde{t})
\end{equation}
where $W(\tilde{t} - t)$ captures the operational energy loss until time $\tilde{t}$, and $\hat{\boldsymbol{z}}$ is P's new position at time $\hat{t} = \tilde{t} + \hat{\tau}$. $\hat{\tau}$ denotes the time it takes for P to first satisfy $|\boldsymbol{z}_E(s) - \boldsymbol{z}_P(s)| = \overline{D}$ flying with velocity $\hat{\boldsymbol{v}}(s) = F_P \left( \frac{\boldsymbol{z}_E(s) - \boldsymbol{z}_P(s)}{|\boldsymbol{z}_E(s) - \boldsymbol{z}_P(s)|}\right)$, and $p$ is computed according to Equation \eqref{eqn:prob_of_no_spook}. 

\section{Numerical Implementation}\label{section:numerics}
Our numerical method is based on a first order semi-Lagrangian (SL) discretization \cite{falcone2013semi} of \eqref{eqn:hjb_stalk_stationary} and \eqref{eqn:hjb_stalk} on a Cartesian grid $Z$ over the domain $\overline{\Omega} := [0, \bar{d}] \times [0, \bar{d}] \times [0, T_{\star}]$. $\bar{d}$ is selected to ensure that $\mathcal{D}(t) \subseteq [0, \bar{d}] \times [0, \bar{d}]$ for all $t \in [0, T_{\star}]$. $Z$ consists of $N_d + 1$, $N_d + 1$, and $N_t + 1$ nodes in the $x, y,$ and $t$ dimensions respectively. For fixed $N_d$,
 the grid spacings are given by
\begin{equation}\label{eqn:grid_spacings}
	\Delta x =  \Delta y = \frac{\bar{d}}{N_d}, \hspace{0.2cm} \Delta t \leq \frac{\sqrt{\Delta x^2 + \Delta y^2}}{F_P}.
\end{equation}
Thus, P's position and time at node $(i,j,k)$ is $(\boldsymbol{z}_{i,j}, t_k) = (i\Delta x, j \Delta y, k \Delta t)$ for $i, j = 0, \dots, N_d$ and $k = 0, \dots, N_t$.

\textit{Stationary HJB Solve:}
Inspired by the approach outlined in \cite{andrews2014deterministic}, we compute the value function at $\boldsymbol{z}_{i,j}$ using a five-point computational stencil constructed from $\boldsymbol{z}_{i,j}$'s four nearest neighbors ($\boldsymbol{z}_{i+1,j}, \boldsymbol{z}_{i,j+1}, \boldsymbol{z}_{i-1,j}, \boldsymbol{z}_{i,j-1}$). Each of the four simplexes in the stencil is a right triangle, and we label them counterclockwise, starting with the northeastern simplex $s^1$ which is characterized by the vertices $ \{\boldsymbol{z}_{i,j}, \boldsymbol{z}_{i+1,j}, \boldsymbol{z}_{i,j+1}\}$.  
Focusing on $s_1$, for all $i, j$, the SL discretization of  \eqref{eqn:hjb_stalk_stationary} satisfies
\begin{multline}
    W^{1}_{i,j} =
\min_{\boldsymbol{\xi} \in \Xi_2} \{K(\boldsymbol{v}^{\boldsymbol{\xi}})\tau(\boldsymbol{\xi}) + pW(\tilde{\boldsymbol{z}}^{\boldsymbol{\xi}}_{i,j}) + \\
 (1-p) \delta(\tilde{\boldsymbol{z}}^{\boldsymbol{\xi}}_{i,j}, T_{\star})\}, \hspace{0.1cm} \boldsymbol{z}_{i,j} \in Z \cap \mathcal{D}(T_{\star}) 
\end{multline}
where $\Xi_2 := \{\boldsymbol{\xi} = (\xi_1, \xi_2) \hspace{0.1cm} | \hspace{0.1cm} \xi_1 + \xi_2 = 1, \hspace{0.1cm} \xi_1, \xi_2 \ge 0 \}$, $\tilde{\boldsymbol{z}}^{\boldsymbol{\xi}}_{i,j} = \xi_1 \boldsymbol{z}_{i+1,j} + \xi_2 \boldsymbol{z}_{i,j+1}$, $\tau(\boldsymbol{\xi}) = \frac{|\tilde{\boldsymbol{z}}^{\boldsymbol{\xi}}_{i,j} - \boldsymbol{z}_{i,j}|}{F_P}$, $\boldsymbol{v}^{\boldsymbol{\xi}}$ is the velocity in the direction of $\tilde{\boldsymbol{z}}^{\boldsymbol{\xi}}_{i,j}$, and $p$ is the probability of E not attempting escape during the $\tau(\boldsymbol{\xi})$ second period.  $W(\tilde{\boldsymbol{z}}^{\boldsymbol{\xi}}_{i,j}) = \xi_1W(\boldsymbol{z}_{i+1,j}) + \xi_2W(\boldsymbol{z}_{i, j+1})$, and any simplex vertices falling within $Z \cap \mathcal{B}(t)$ are replaced with the nearest point visible to E. The value function at $\boldsymbol{z}_{i,j}$ is
\begin{equation}
	W_{i,j} = \min\{W^1_{i,j}, W^2_{i,j}, W^3_{i,j}, W^4_{i,j}\},
\end{equation}
and the boundary condition on $\mathcal{T}$ is $W_{i,j} = \delta(\boldsymbol{z}_{i,j},T_{\star})$. 

\textit{Time-Dependent HJB:} Equation \eqref{eqn:hjb_stalk} is treated similarly, but with slight modifications to accommodate the time dependence. We assume that P at $\boldsymbol{z}_{i,j}$ applies his choice of control $\boldsymbol{v}$ for $\tau$ seconds before arriving at $\tilde{\boldsymbol{z}}^{\boldsymbol{v}}_{i,j} \approx \boldsymbol{z}_{i,j} + \tau \boldsymbol{v}$. The SL discretization of \eqref{eqn:hjb_stalk} for all $i, j,$ and $k < N_t$ satisfies 
\begin{multline}\label{eqn:sl_discretization}
    U^k_{i,j} =  \min_{\boldsymbol{v} \in A_F}\{\tau K(\boldsymbol{v}) + pU(\tilde{\boldsymbol{z}}^{\boldsymbol{v}}_{i,j}, t + \tau) \\
    + (1-p)\tilde{\delta}(\boldsymbol{z}^{\boldsymbol{v}}_{i,j}, t+\tau)\}, \hspace{0.1cm} \boldsymbol{z}_{i,j} \in Z \cap \mathcal{D}(t_k) 
\end{multline}
where $\tau = \Delta t$, and $p \approx \tau \lambda(\boldsymbol{z}_{i,j}, t_k)$. We solve backwards in time from $U^{N_t}_{i,j} = W_{i,j}$, and we maintain the boundary condition $U^k_{i,j} = \delta(\boldsymbol{z}_{i,j}, T_1)$, $\forall \boldsymbol{z}_{i,j} \in Z \cap \mathcal{T}$. On portions of the domain where $|\tilde{\boldsymbol{z}}^{\boldsymbol{v}}_{i,j} - \boldsymbol{z}_E(t)| \leq D$, we recover $U(\tilde{\boldsymbol{z}}^{\boldsymbol{v}}_{i,j}, t + \tau)$ via bilinear interpolation in $x$ and $y$. 
At points $\boldsymbol{z}_{i,j} \in Z \cap \mathcal{D}(t_k)$ which are outside of E's visual range at time $t_k$, we compute $U^k_{i,j}$ using the approximation formula in \eqref{eqn:visual_distance_cost} with trilinear interpolation in $x, y$, and $t$ to determine $U(\hat{\boldsymbol{z}}, \hat{t})$.

\textit{Modified Hit-and-Run Algorithm for $\boldsymbol{v}_*$:}
While the optimal velocity $\boldsymbol{v}_*$ is known a priori at nodes $Z \cap \mathcal{T}$ and at nodes satisfying $|\boldsymbol{z}_{i,j} - \boldsymbol{z}_E(t)| \ge D$, we utilize a derivative-free Hit-and-Run optimization algorithm modeled after the one proposed in \cite{zabinsky1993improving} to determine $\boldsymbol{v}_*$ at nodes which fall in the annulus $\epsilon < |\boldsymbol{z}_{i,j} - \boldsymbol{z}_E(t)| < D$. We first perform a coarse grid search by discretizing $A_F$ along $a+1$ angular directions $\theta_{\ell} = \frac{2\pi \ell}{a}\in [0, 2\pi)$, $\ell = 0, \dots, a$ and $b+1$ speed values, $s_n = \frac{F_p n}{b} \in [0, F_P]$, $n = 0, \dots, b$. For each velocity $\boldsymbol{v}_{\ell n} = [s_n\cos(\theta_{\ell}), s_n\sin(\theta_{\ell})]^T$, we obtain a value function candidate $V_{\ell n}(\boldsymbol{z}_{i,j}, t_k)$ computed according to the formula in \eqref{eqn:sl_discretization}. We then execute a Hit-and-Run algorithm with the initial objective function value $V^* = \min\left\{V_{\ell n}\right\}$. The algorithm uniformly samples a search direction $\tilde{\theta} \in [0, 2\pi)$ and step size $\tilde{\gamma} \in [0,1]$ to generate a new velocity $\tilde{\boldsymbol{v}} = [s_n\cos(\theta_{\ell}) + \tilde{\gamma}\cos(\tilde{\theta}), s_n\sin(\theta_{\ell}) + \tilde{\gamma}\sin(\tilde{\theta})]^T$, and computes the corresponding value function candidate $\tilde{V}$. If $\tilde{V} < V^*$, we set $V^* = \tilde{V}$, and repeat this procedure until improvements between successive accepted optimum candidates falls below a threshold or the search exceeds a max number of iterations. 

\textit{Optimal Trajectory Tracing \& MC Optimality:} After computing $U$ and $\boldsymbol{v}_*$ at each $(\boldsymbol{z}_{i,j}, t_k)$, 
we reconstruct an approximation to P's full optimal trajectory until E's capture starting from a specified $(x, y)$ at $t = 0$. We approximate $\boldsymbol{v}_*(\boldsymbol{z}, t)$ via bilinear interpolation in $x$ and $y$ and integrate P's stalk phase dynamics in steps of $\Delta t$ up until $t = T_s$. Since it is impossible to recognize perfect MC due to computational and physiological limitations, we determine the portions of P's optimal trajectory which are \emph{approximately MC} (A-MC) by calculating $\theta_{\#}$ using the formula given in Equation \eqref{eqn:theta_sharp} at each point along the stalk phase trajectory. A-MC occurs when $\theta_{\#} \leq 0.5^{\circ}$, which is below the resolvable angle limit for many compound eye visual systems. From $T_s$ onward, we integrate P's direct chase phase dynamics along the analytically known optimal direction of motion $\boldsymbol{a}_*(t)$ until P captures E at $t = T_2$. 

\section{Numerical Experiments}\label{section:experiments}
In all examples\footnote{Source code, additional figures, and movies are available at \url{https://github.com/eikonal-equation/Motion_Camo_Optimality}.}, E's stalk phase trajectory is given by $\boldsymbol{f}(t) = [2 + 0.05t - \cos(t), 1 + 0.05 t^2 + \sin(t)]^T$, and $T_{\star} = 2s$. The surroundings are homogeneous except for a boulder at $\boldsymbol{z}_{\#} = (1.5, 0.6)$.  Following experimental observations in \cite{collett1978hoverflies}, we set $D = 0.75 m$, $\bar{D} = 0.15 m$, and $G_E = 10 m/s$.  We take $\gamma = 0.025m$, $\epsilon = 0.05m$, $F_P = 4 m/s$, $G_P = 20 m/s$, and $W$ is set to $3 J / s$. We impose a $201 \times 201 \times 801$ grid over the domain $\overline{\Omega} = [0, 4] \times [0, 4] \times [0, 2]$, and the spacings are computed according to the formulas in \eqref{eqn:grid_spacings}. In the Modified Hit-and-Run for $\boldsymbol{v}_*$, $a = 30$, $b = 15$, and it terminates when the error between successive accepted optimum candidates is less than $1\times 10^{-4}$ or it has exceeded $100$ iterations.

\begin{figure*}[h!]
$
\arraycolsep=1pt\def\arraystretch{0.1}
\begin{array}{cccc}
\includegraphics[width=0.25\textwidth]{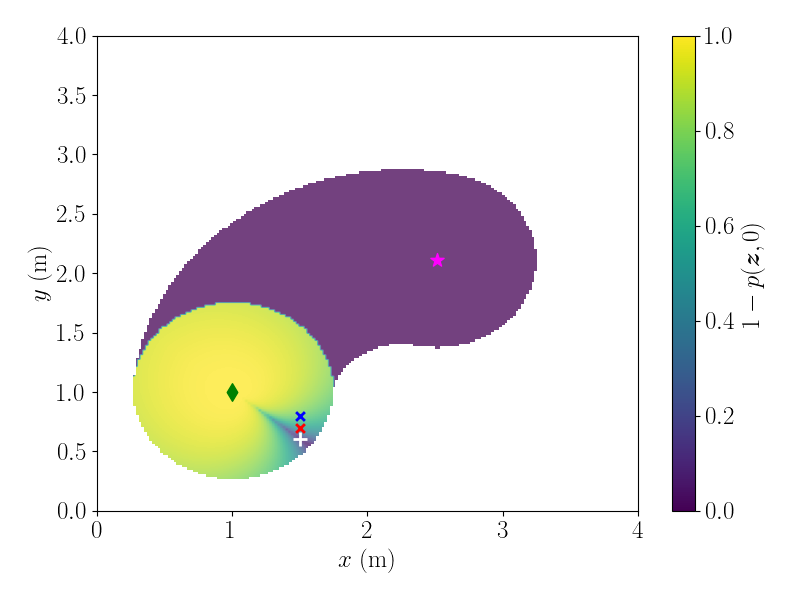} &
\includegraphics[width=0.25\textwidth]{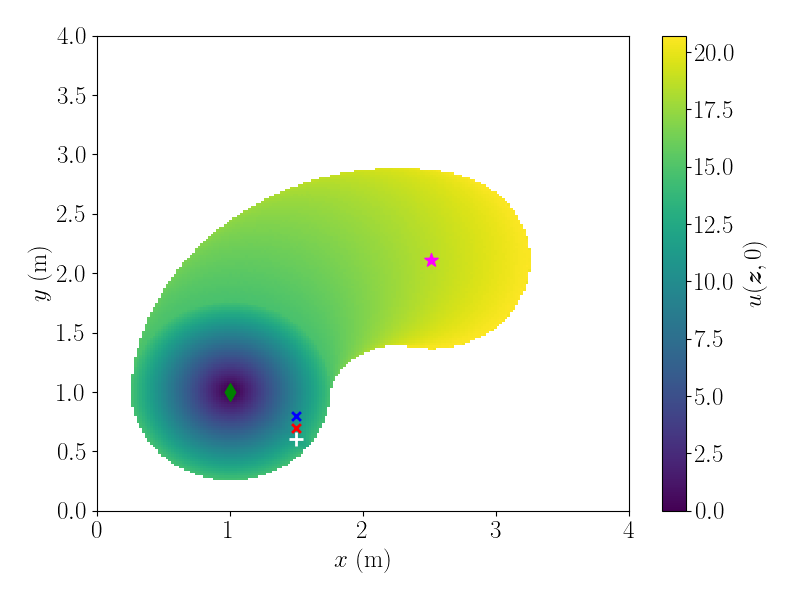} &
\includegraphics[width=0.25\textwidth]{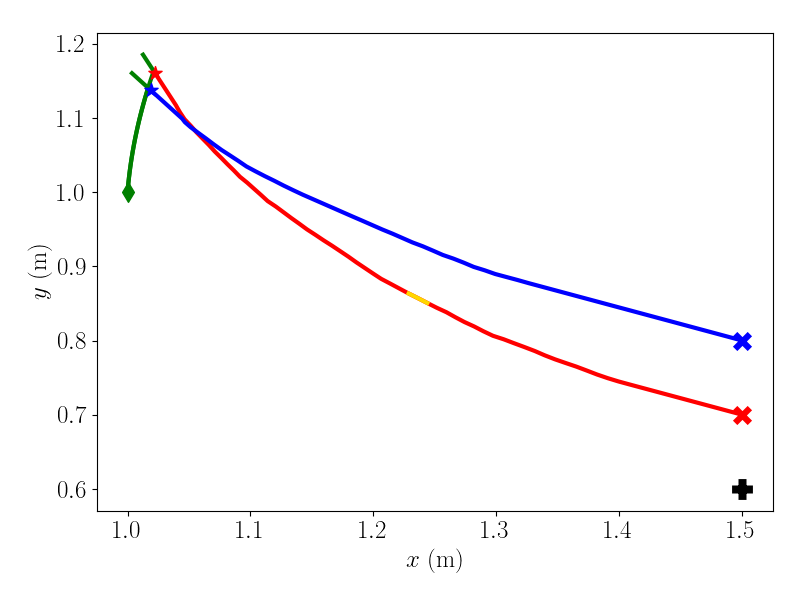} &
\includegraphics[width=0.25\textwidth]{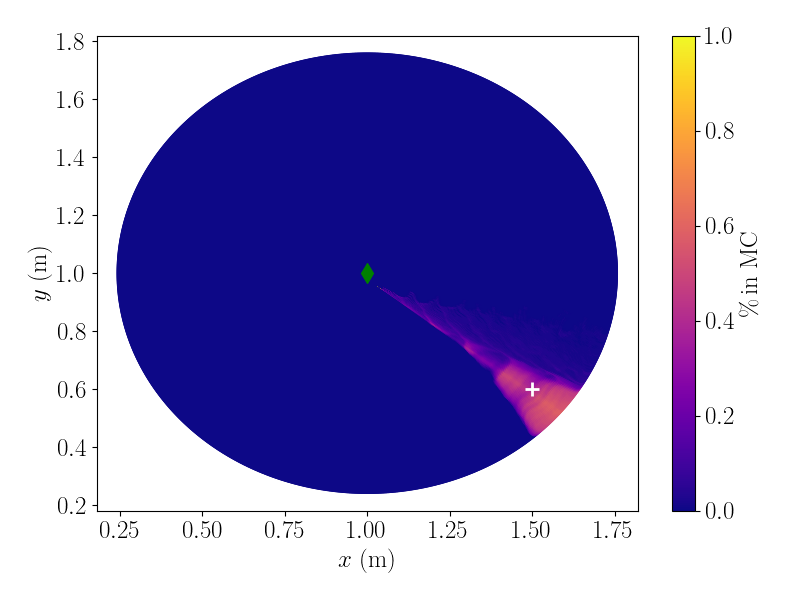}\\[-3pt]
\mbox{\footnotesize (A)} & \mbox{\footnotesize(B)} & \mbox{\footnotesize(C)} & \mbox{\footnotesize(D)}
\end{array}
$
\caption{Example 1 pointwise detection probability at $t = 0$ (A), value function at $t = 0$ (B), optimal trajectories starting from $(1.5, 0.7)$ and $(1.5, 0.8)$ (C), and $\%$ of flight time spent in approximate-MC starting within $D \hspace{0.1cm} m$ of E at $t=0$ (D). Given E's low tolerance and low visual acuity, adjusting to an approximately-MC trajectory is not worthwhile. The approximately-MC portion of P's trajectory starting from $(1.4, 0.7)$ is merely coincidental. \label{fig:ex1}}
\end{figure*}

\begin{figure*}[h!]
$
\arraycolsep=1pt\def\arraystretch{0.1}
\begin{array}{cccc}
\includegraphics[width=0.25\textwidth]{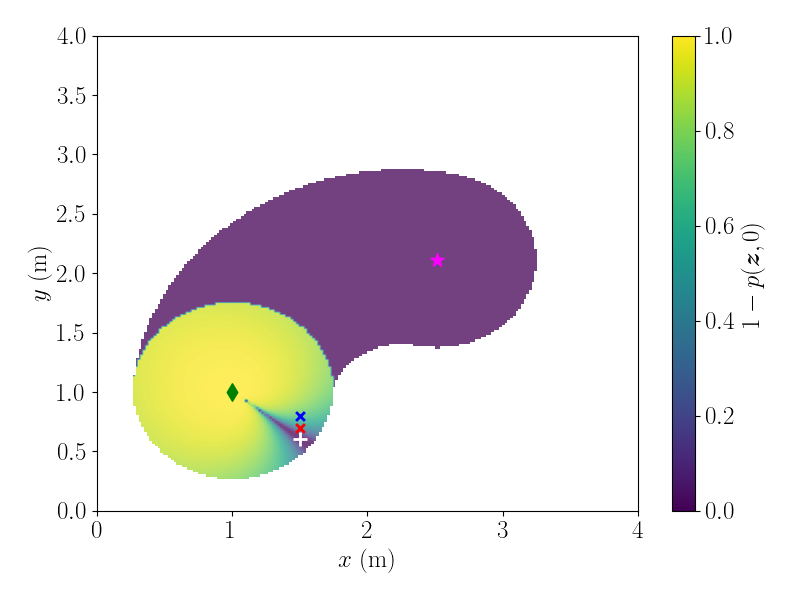} &
\includegraphics[width=0.25\textwidth]{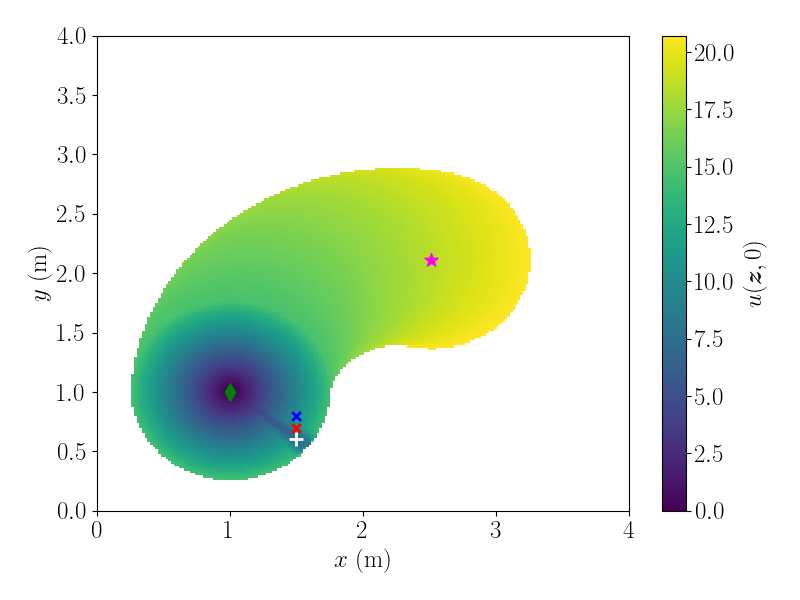} &
\includegraphics[width=0.25\textwidth]{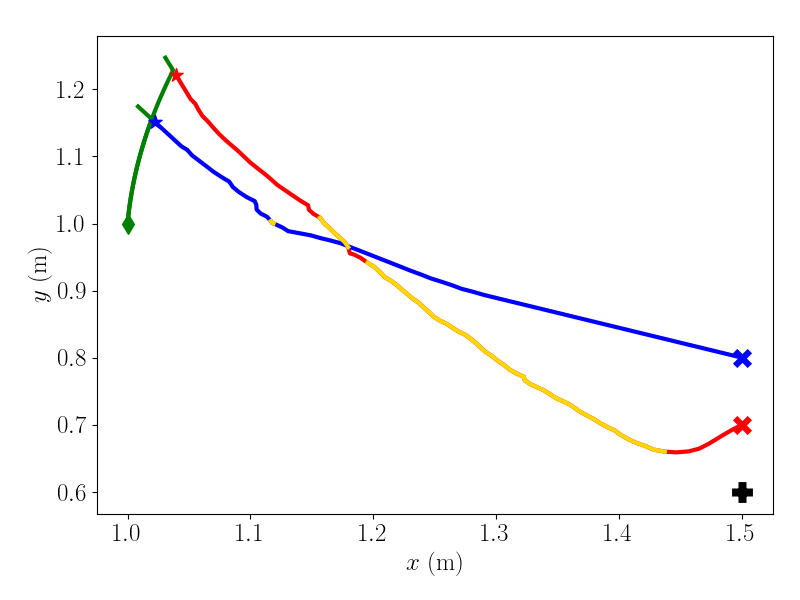} &
\includegraphics[width=0.25\textwidth]{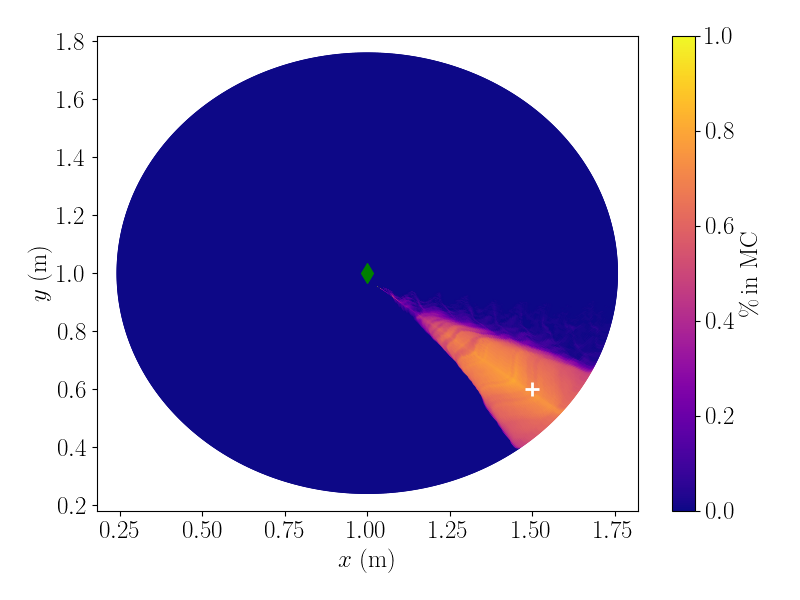}\\[-3pt]
\mbox{\footnotesize (A)} & \mbox{\footnotesize(B)} & \mbox{\footnotesize(C)} & \mbox{\footnotesize(D)}
\end{array}
$
\caption{Example 2 pointwise detection probability at $t = 0$ (A), value function at $t = 0$ (B), optimal trajectories starting from $(1.5, 0.7)$ and $(1.5, 0.8)$ (C), and $\%$ of flight time spent in approximate-MC starting within $D \hspace{0.1cm} m$ of E at $t=0$ (D). P spends $64\%$ and $6.1\%$ of his stalk phase flight
time in approximate-MC starting from (1.5, 0.7) and (1.5, 0.8)
respectively.  Despite low tolerance, E's strong visual acuity makes the energy required to adjust to an A-MC trajectory worthwhile for P starting near and within the small sector where E's escape attempt probability is $\approx 0$. \label{fig:ex2}}
\end{figure*}

\begin{figure*}[h!]
$
\arraycolsep=1pt\def\arraystretch{0.1}
\begin{array}{cccc}
\includegraphics[width=0.25\textwidth]{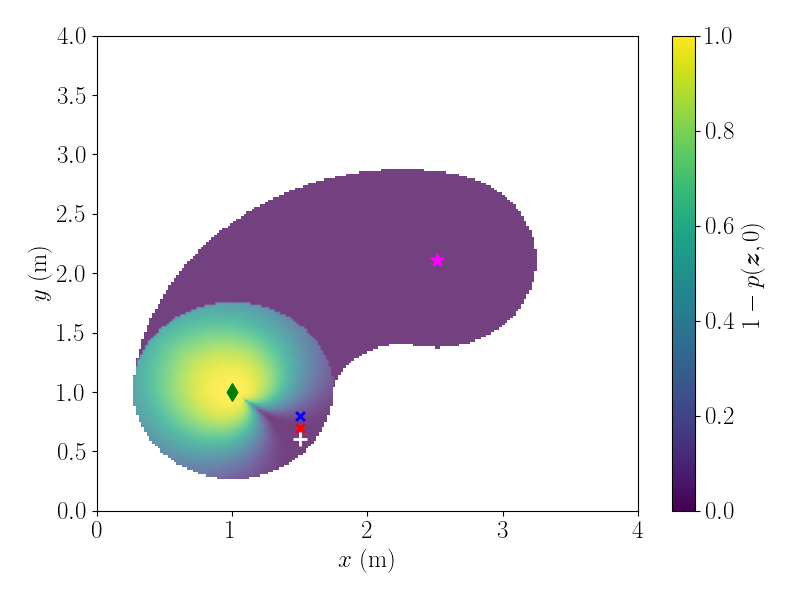} &
\includegraphics[width=0.25\textwidth]{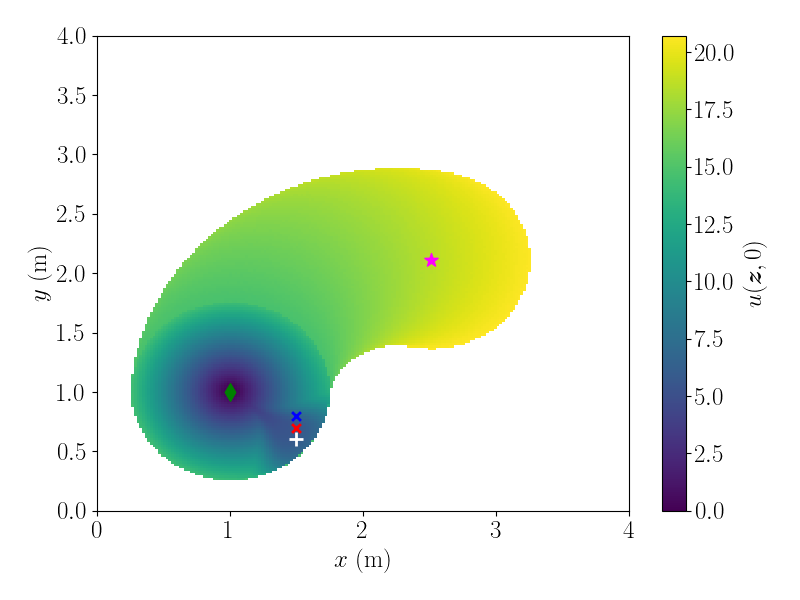} &
\includegraphics[width=0.25\textwidth]{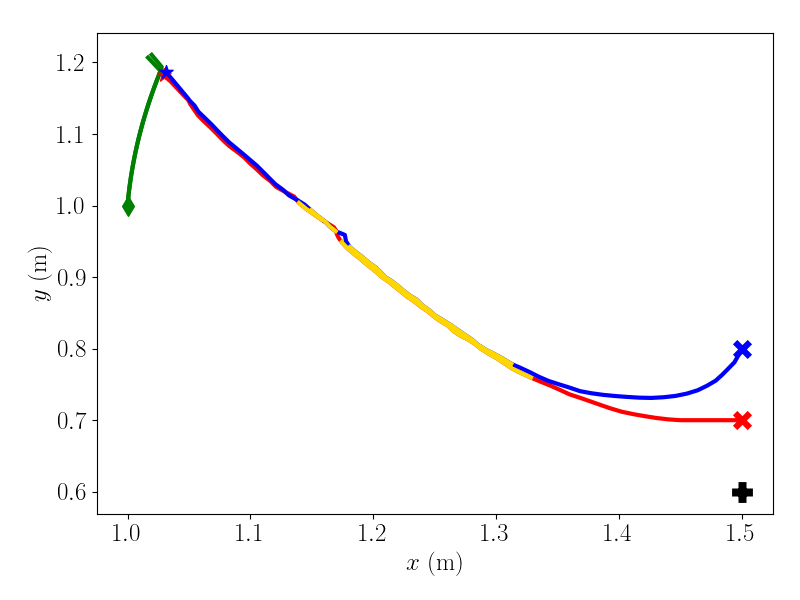} &
\includegraphics[width=0.25\textwidth]{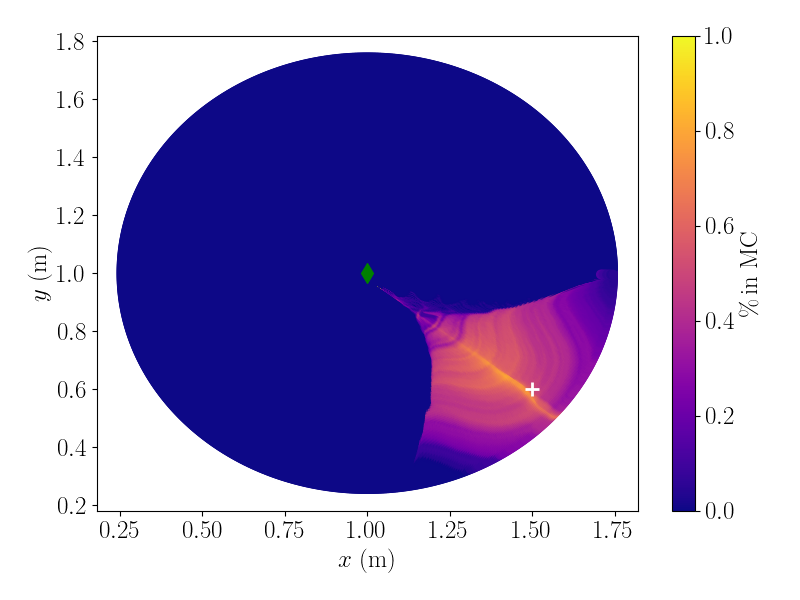}\\[-3pt]
\mbox{\footnotesize (A)} & \mbox{\footnotesize(B)} & \mbox{\footnotesize(C)} & \mbox{\footnotesize(D)}
\end{array}
$
\caption{Example 3 pointwise detection probability at $t = 0$ (A), value function at $t = 0$ (B), optimal trajectories starting from $(1.5, 0.7)$ and $(1.5, 0.8)$ (C), and $\%$ of flight time spent in approximate-MC starting within $D \hspace{0.1cm} m$ of E at $t=0$ (D). P spends $49\%$ and $47\%$ of his stalk phase flight
time in A-MC starting from (1.5, 0.7) and (1.5, 0.8)
respectively. Regardless of E's low visual acuity, higher tolerance for P significantly lowers the chance of premature switching to the chase phase from most locations within E's visual range, and MC tactics are worthwhile along more trajectories (e.g., from $(1.5, 0.8)$).  \label{fig:ex3}}
\end{figure*}

\begin{figure*}[h!]
$
\arraycolsep=1pt\def\arraystretch{0.1}
\begin{array}{cccc}
\includegraphics[width=0.25\textwidth]{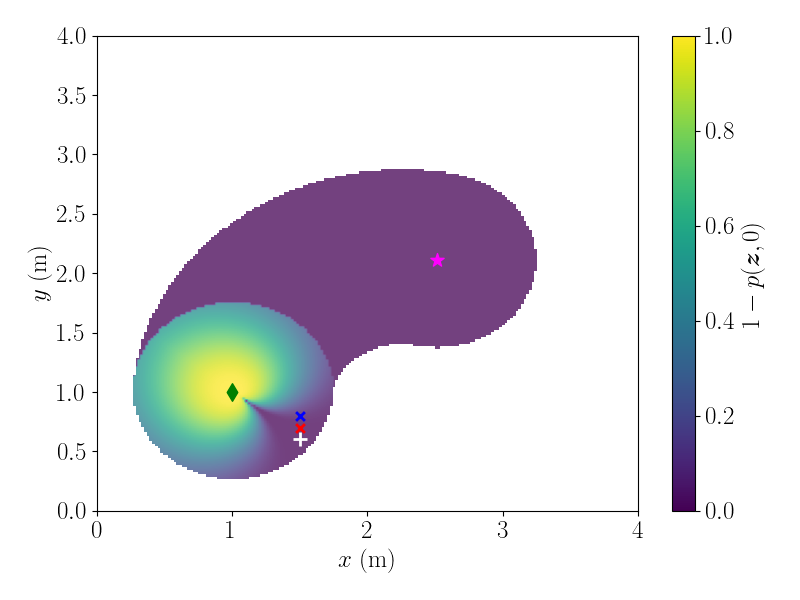} &
\includegraphics[width=0.25\textwidth]{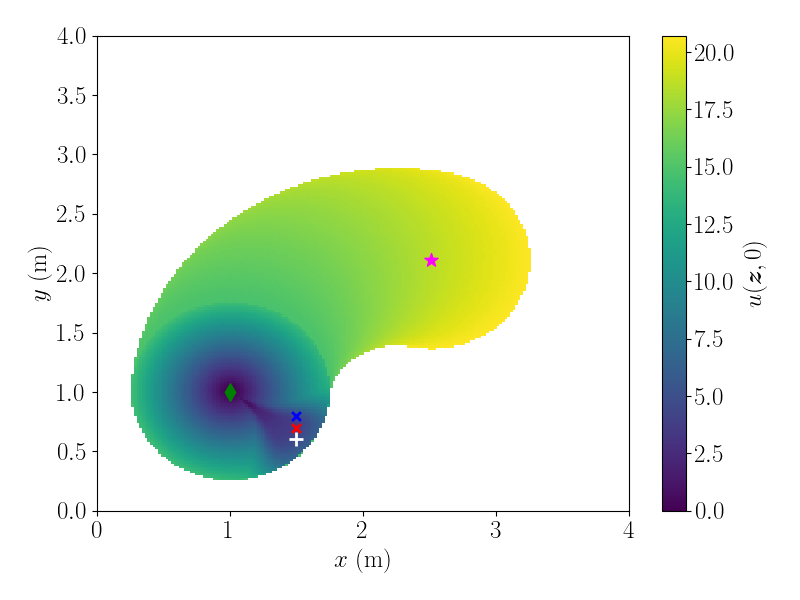} &
\includegraphics[width=0.25\textwidth]{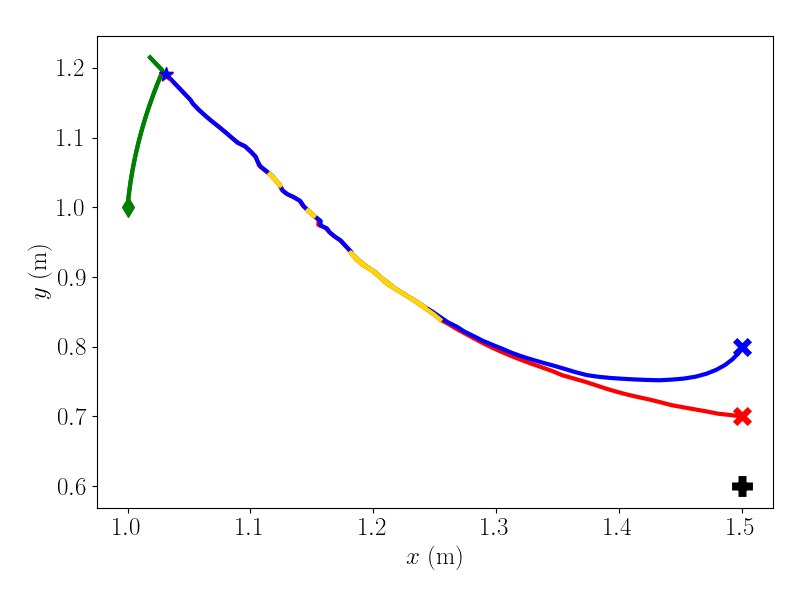} &
\includegraphics[width=0.25\textwidth]{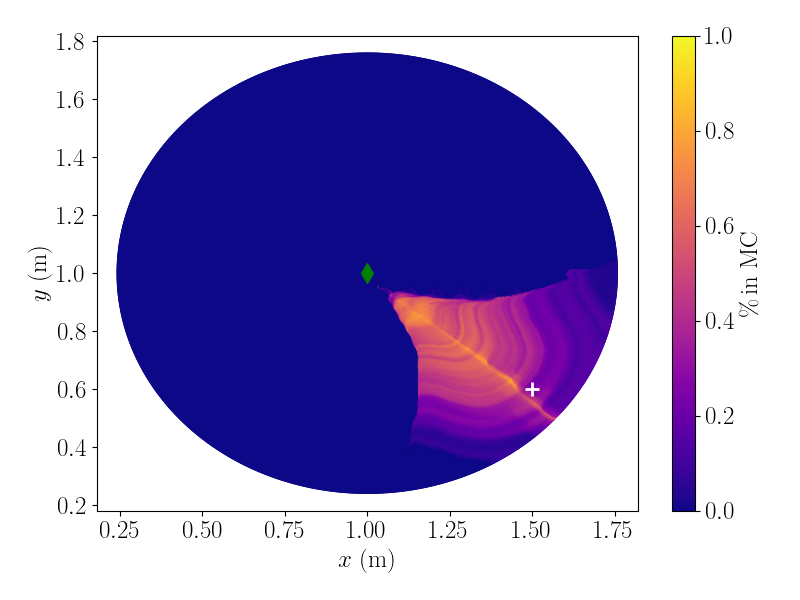}\\[-3pt]
\mbox{\footnotesize (A)} & \mbox{\footnotesize(B)} & \mbox{\footnotesize(C)} & \mbox{\footnotesize(D)}
\end{array}
$
\caption{Example 4 pointwise detection probability at $t = 0$ (A), value function at $t = 0$ (B), optimal trajectories starting from $(1.5, 0.7)$ and $(1.5, 0.8)$ (C), and $\%$ of flight time spent in approximate-MC starting within $D \hspace{0.1cm} m$ of E at $t=0$ (D). From both starting points, P spends $38\%$
of his stalk phase flight time in A-MC. Although E's visual acuity is high, E's higher tolerance for P makes engaging in MC for a significant portion of the flight less worthwhile for trajectories starting farther away from E at $t=0$, but more worthwhile for trajectories starting closer to E at $t = 0$. \label{fig:ex4}}
\end{figure*}

Our experiments illustrate how E's visual acuity and tolerance affects the optimality of MC tactics during the stalk phase by numerically solving the multistage optimal control problem for various values of $B$ and $C$ in \eqref{eqn:lambda}. Parameter values and the resulting percentage of trajectories starting within E's visual range for which approximate-MC is optimal for over $50\%$ of the flight time are summarized for the four examples in Table \ref{table:mc_stats}. Figure \ref{fig:ex1} provides a grid of plots for each example. In all plots, $\boldsymbol{z}_{\#}$ is indicated by a ``$+$", and E's starting position at $t = 0$ is marked by a green diamond.  The pointwise detection probability at $t=0$ is shown in Figures \ref{fig:ex1} - \ref{fig:ex4} (A). The value function at $t = 0$ is shown in Figures \ref{fig:ex1} - \ref{fig:ex4} (B). Figures \ref{fig:ex1} - \ref{fig:ex4} (C) display optimal trajectories starting from the red ``X" at $(1.5, 0.7)$ and the blue ``X" at $(1.5, 0.8)$ without premature switching to direct chase. E's full trajectory until capture is shown in green, the end of P's stalk phase trajectory is marked by a ``$\star$", and approximate-MC portions of P's optimal trajectory are highlighted in gold. Figures \ref{fig:ex1} - \ref{fig:ex4} (D) illustrate the percentage of P's stalk phase flight time spent in A-MC starting from $200,000$ unique points visible to E at $t = 0$.

\begin{table}
\caption{Experimental Parameters \& Summary Statistics}
\label{table:mc_stats}
\begin{tabular}{ |p{1.5cm}||p{1.5cm}|p{1.5cm}|p{2cm}|  }
 \hline
 \multicolumn{4}{|c|}{$\%$ of Trajectories for which A-MC is Optimal More Than $50\%$ of Flight} \\
 \hline
Example & $B$ & $C$ & \% of Trajectories \\
 \hline
1  & $0.05$    & $0.4$ & $0.18\%$\\
 2&   $5 \times 10^{-5}$ & $0.4$   & $3.4\%$\\
 3 & $0.05$ & $4$ &  $2.7\%$ \\
4    & $5 \times 10^{-5}$ & $4$ &  $4.3\%$\\
 \hline
\end{tabular}
\end{table}


\section{Conclusions}\label{section:conclusions}
We presented a continuous-time dynamic programming framework to determine when it is optimal for an energy-optimizing pursuer to utilize motion camouflage amidst uncertainty in the evader's escape response. We illustrated our setup for the biological example of a male hover fly pursuing a female hover fly, and we showed how varying visual acuity and tolerance to P's presence in the female's visual range can affect the stalk-phase optimality of MC strategies. 

Despite our focus on hover fly interactions, we emphasize that this approach to analyzing motion camouflaging behaviors is applicable to a broader class of pursuer-evader interactions. One can readily modify the control, cost function, dynamics, and $\lambda$ as appropriate to apply the framework to other biological or vehicular systems. Possible extensions include adding more food sources to introduce uncertainty in E's intended destination and analyzing how properties of E's trajectory affect MC-optimality. We also hope that our approach may serve as a building block for optimization under uncertainty in more complex scenarios such as optimal pursuit in multi-pursuer / multi-evader systems \cite{zhou2024optimal}, motion camouflage in mutual pursuit systems \cite{mischiati2012dynamics}, and guiding bio-inspired robots in three spatial dimensions \cite{savkin2020bioinspired}. 

\section*{Acknowledgment}
The author thanks Alex Vladimirsky for helpful discussions and the reviewers for useful feedback.

\vspace{-0.2cm}
\bibliographystyle{IEEEtran}
\bibliography{bibl}

\end{document}